\documentclass[12pt]{amsart}
\usepackage{amssymb,amsmath,latexsym,amscd,amsfonts}
\usepackage{pgf,tikz,pgfplots}
\pgfplotsset{compat=1.14}
\usepackage{mathrsfs}
\usetikzlibrary{arrows}
\usetikzlibrary[patterns]
\usepackage{graphics}
\usepackage{epsfig}
\usepackage{psfrag}
\usepackage{wasysym}

\usepackage[pagewise]{lineno}%\linenumbers

\def\ignore #1 {}

\newtheorem{thm}{Theorem}
\newtheorem{lem}[thm]{Lemma}

\newtheorem{prop}[thm]{Proposition}
\newtheorem{cor}[thm]{Corollary}

\theoremstyle{definition}
\newtheorem*{remark}{Remark}

\newtheorem*{example*}{Example}

\def\hpic #1 #2 {\mbox{$\begin{array}[c]{l} \epsfig{file=#1,height=#2} \end{arr\
ay}$}}
\def\vpic #1 #2 {\mbox{$\begin{array}[c]{l} \epsfig{file=#1,width=#2} \end{arra\
y}$}}

\def\Q{\mbox{{\bf Q}}}

\def\C{\mbox{{\bf C}}}

\def\Z{\mbox{{$\mathbb{Z}$}}}

\DeclareMathOperator{\im}{Im}

\DeclareMathOperator{\Area}{Area}
\DeclareMathOperator{\Vxs}{Vertices}

\DeclareMathOperator{\Trap}{Trap}
\DeclareMathOperator{\Frac}{Frac}

\newcommand{\s}{\sigma}

\def \sb #1 {\overline{\s}_{#1}}

\def\PP{\mathbf{p}}
\def\QQ{\mathbf{q}}
\def\RR{\mathbf{r}}
\def\SS{\mathbf{s}}
\def\thering{\mathcal{R}}

\begin{document}
%%%%%%%%%%%%%% TITLE STUFF %%%%%%%%%%%%%%%%%%%%
\title{Integrality relations for polygonal dissections}
\author{Aaron Abrams}
\address{Washington and Lee University}
\email{\tt abramsa@wlu.edu}
\author{Jamie Pommersheim}
\address{Reed College}
\email{jamie@reed.edu}
%
%\keywords{}

\maketitle

\begin{abstract} 
  Given a parallelogram dissected into triangles, the area of any one of the triangles of the dissection is integral over the ring generated by the areas of the other triangles. Given a trapezoid dissected into triangles, the area of any triangle determined by either diagonal of the trapezoid is integral over the ring generated by the areas of the triangles in the dissection.  In both cases, the integrality relations are invariant under deformation of the dissection.
  
  The trapezoid theorem implies and provides a new context for Monsky's Equidissection Theorem that a square cannot be dissected into an odd number of triangles of equal area. 
  A corollary of these results is that the area polynomials for parallelograms introduced in \cite{triangles1} and studied further in \cite{chapter3,chapter2} have all leading coefficients equal to $\pm 1$.  
  
\end{abstract}

\bigskip

\section{Introduction}

We establish several new results about the geometry of dissections of certain Euclidean plane polygons.
A \emph{dissection} of such a polygon $T$ into triangles is a collection of triangles in the plane whose union is $T$ and whose interiors are disjoint.

\begin{thm}\label{thm1}
Let $T$ be a trapezoid in the Euclidean plane with vertices $\PP,\QQ,\RR,\SS$, in counterclockwise order. Suppose that $T$ is  dissected into $n$ triangles of areas $a_1, \dots, a_n$.  Then the area of the triangle $\PP\QQ\SS$ is integral over $\Z[a_1,\dots, a_n]$.
\end{thm}

\begin{thm}\label{thm2}
Let $T$ be a parallelogram in the Euclidean plane with a dissection into $n$ triangles of areas $a_1, \dots, a_n$. Then $a_n$ is integral over $\Z[a_1,\dots, a_{n-1}]$. 
\end{thm}

Theorem \ref{thm1} immediately implies Monsky's theorem \cite{monsky} that a parallelogram cannot be dissected into an odd number of triangles of equal area, since $1/2$ is not integral over $\Z[1/n]$ when $n$ is odd. Thus Theorem \ref{thm1} generalizes and provides a new context for Monsky's Theorem. However, this cannot be considered a new proof of Monsky's Theorem, since our proof proceeds along the same lines as the original, using valuations to 3-color points of a certain affine plane and appealing to Sperner's Lemma.  See \cite{monsky,monskyjepsen}.

We also show that in a certain sense, the integrality relations arising in these theorems are invariant under deformation; that is, the integrality relations actually hold for the quadratic polynomials that express the areas of the triangles, and not just for the numerical areas $a_i$.  See Theorem 1+ and Theorem \ref{monicity} below.

Theorem \ref{thm2} goes hand in hand with a result about the {\it area polynomial} $p_T$ that was introduced in \cite{triangles1} and further studied in \cite{chapter3} and \cite{chapter2}.  For any combinatorial triangulation $T$ of a quadrilateral, there is a unique (up to sign) nonzero homogeneous irreducible integer polynomial $p_T$ with one variable $A_i$ for each triangle such that $p_T(a_1, \dots, a_n)=0$ whenever $T$ is drawn in the plane with a parallelogram boundary and triangles of areas $a_1, \dots, a_n$. Here by \emph{combinatorial triangulation} we mean a simplicial complex homeomorphic to a disk, with four vertices on the boundary. (The connection with dissections is that every dissection of a planar trapezoid can be viewed as the image of a combinatorial triangulation under a piecewise linear map to the plane which may collapse some triangles; see e.g.~\cite[Propositions 2 and 5]{chapter3}.)
The mod 2 structure of $p_T$ is completely specified by \cite[Theorem 9.1]{chapter3}, which implies in particular that the coefficients of the leading terms are odd integers.  Further, in \cite[Theorem 6.2]{chapter2} it is shown that these leading terms must all be equal up to sign.  %Theorem \ref{thm2} in the present paper allows us to strengthen this.

\begin{thm}\label{monicity}
For any combinatorial triangulation $T$, the area polynomial $p_T$ is monic.  That is, for any $i$ the coefficient of $A_i^d$ is $\pm 1$, where $d=\deg p_T$.
\end{thm}

This is a special case of the positivity conjecture from \cite[Conjecture 4]{chapter3}.

\begin{remark}
Monsky's equidissection theorem, as well as our Theorems 1 and 2, apply to arbitrary dissections, whereas the combinatorial triangulations of Theorem 3 are by definition simplicial complexes.   
It is easy to see that Theorem 3 also holds for any dissection that has an area polynomial (the ``hyper" case in the language of \cite[Defintion 26]{chapter3}). However it is not known whether every dissection of a parallelogram has this property; this question is discussed in Section 8 of \cite{chapter3}.
\end{remark}

We also note that integrality conditions have previously appeared in theorems about equidissections of trapezoids.  For example \cite[Theorem 1.1]{monskyjepsen} (see also \cite{kasimatisstein}) gives a necessary condition for the existence of an equidissection of a trapezoid of a given shape into a given number of triangles.  Theorem \ref{thm1} strengthens that result.

The theorems are proved by combining ideas originally due to Monsky \cite{monsky} with some technical machinery developed in \cite{triangles1,chapter3,chapter2}.  Some familiarity with those works may be helpful for the reader; in order to focus on the results, we have not attempted to make the arguments here entirely self-contained.

\bigskip

{\bf Acknowledgements.}
The authors thank Paul Monsky for numerous insightful communications and inspirations.  In addition the first author gratefully acknowledges the support of the MTA Distinguished Guest Scientist Fellowship Programme 2022.
The second author thanks the Fulbright U.S. Scholar Program for their support. We also thank Dezs\H{o} Mikl\'os and Andr\'as Stipsicz for their roles in making this work possible.

\section{Integrality for Trapezoids}

In this section, we prove Theorem \ref{thm1} by establishing an integrality relation for the triangle $\PP\QQ\SS$ of a dissected trapezoid.  In fact we prove a stronger version of this theorem (Theorem 1+) that allows deformations of the trapezoid. 

Let $T$ be an combinatorial triangulation of a quadrilateral $\PP\QQ\RR\SS$.  For each vertex $v$ other than $\RR$, we introduce two variables $x_v$ and $y_v$.  We treat $v=\RR$ differently so that our ring will reflect the geometric condition that $\PP\QQ\RR\SS$ be a trapezoid rather than an arbitrary quadrilateral. For this final vertex, we introduce a variable $t$ which represents the ratio of the lengths of side $\SS\RR$ to side $\PP\QQ$. Thus we work in the polynomial ring
$$
\thering = \C[\{x_v, y_v | v\in \Vxs(T)\setminus \{\RR\} \}, t].
$$
In $\thering$, we use the abbreviations $x_{\RR} = x_{\SS} + t(x_{\QQ}-x_{\PP})$ and  $y_{\RR} = y_{\SS} + t(y_{\QQ}-y_{\PP})$. In $\thering$, it is natural to consider the variables $x_v$ and $y_v$ as having degree 1, while $t$ has degree $0$.

Orienting the boundary in the direction $\PP\QQ\RR\SS$ endows each triangle $\Delta_i$ of $T$ with an orientation. 
For each $\Delta_i$, we introduce a quadratic polynomial $W_i\in\thering$ which expresses twice the area of the oriented triangle $\Delta_i$.  For convenience, we prefer to work with doubled areas throughout.  This makes little difference, as all the relations we obtain will be homogeneous. We use $W_U\in\thering$ to denote the quadratic polynomial representing twice the area of triangle $\PP\SS\QQ$; this choice of orientation is consistent with the other triangles.  We sometimes abuse language and refer to the $W_i$ and $W_U$ as the areas.

\begin{thm}[Theorem 1+]\label{thm1+}
Let $T$ be a combinatorial triangulation of a quadrilateral $\PP\QQ\RR\SS$ into $n$ triangles.  Let $W_1, \dots, W_n\in\thering$ denote the polynomials expressing the areas of the triangles of $T$, and let $W_U\in\thering$ denote the polynomial expressing the area of the triangle $\PP\SS\QQ$.  Then $W_U$ is integral over $\Z[W_1, \dots, W_n]$.
\end{thm}

\begin{proof}
We use many of the ideas from the proof of Theorem 7.2 (Monsky+) from \cite{chapter3}.  
To show that $W_U$ is integral over $\mathcal{S}=\Z[W_1, \dots, W_n]$, it is enough to show that if $\nu$ is a valuation on the fraction field of $\Z[W_U, W_1, \dots, W_n]$ such that $\nu(W_i)\geq 0$ for all $i$, then $\nu(W_U)\geq 0$ (see, e.g, \cite[5.22]{atiyahmacdonald}). Given such a $\nu$, extend it to the fraction field $\mathcal{F}=\Frac(\thering)$ and, following Monsky \cite{monsky}, use $\nu$ to color each point of $\mathcal{F}\times \mathcal{F}$ one of three colors $A,B,C$ as in the proof of \cite[Theorem 7.2]{chapter3}.

Let $M:\mathcal{F}\times \mathcal{F}\to \mathcal{F}\times \mathcal{F}$ be the unique affine transformation taking $(x_\PP,y_\PP)$ to $(0,0)$, $(x_\QQ,y_\QQ)$ to $(1,0)$, and $(x_\SS,y_\SS)$ to $(0,1)$.
Note that the determinant of $M$ is
$\left|
\begin{matrix}
x_\QQ-x_\PP & x_\SS-x_\PP \\
y_\QQ-y_\PP & y_\SS-y_\PP
\end{matrix}
\right|^{-1}$,
which equals $-W_U^{-1}$.  
 
We now color the vertices of $T$ by using $M$ to pull back the coloring of $\mathcal{F}\times \mathcal{F}$. That is, if $v$ is a vertex of $T$, then we color $v$ with the color of the point $M(x_v,y_v)$. This assigns $\PP, \QQ, \SS$ the colors $C, A, B$, respectively.  As for $\RR$, one sees that $M(x_{\RR}, y_{\RR})= (t,1)$, so $\RR$ has color $A$ or $B$. The boundary of $T$ is thus colored $CAAB$ or $CABB$, and in either case we may apply Sperner's Lemma to conclude that $T$ has an $ABC$ triangle $\Delta_j$. For such a triangle we have $\nu(\Area(M\Delta_j))\leq 0$, which means $\nu(-W_U^{-1}W_j)\leq 0$.  Hence $\nu(W_j)\leq \nu(W_U)$, which implies $\nu(W_U)\geq 0$. 
\end{proof}

We now show that Theorem 1+ implies Theorem \ref{thm1}. 
\begin{proof}
Let $\Delta=\PP\QQ\RR\SS$ be a trapezoid in the plane with a dissection into $n$ triangles of areas $a_1, \dots, a_n$, and let $u$ denote the area of triangle $\PP\SS\QQ$. As in \cite[Propositions 2.6, 3.2]{chapter3}, there exists a combinatorial triangulation $T$ with $m\ge n$ triangles obtained by poofing the dissection, and a drawing $\rho$ of $T$ that has the same set of nondegenerate triangles as the original dissection along with $m-n$ degenerate triangles of area $0$. By Theorem 1+, there is an integral equation $g_T(W_U,W_1,\dots,W_m)=0$, where we may take $g_T$ to be homogeneous in its $m+1$ variables. If $u=0$, then we are done. Otherwise, $\rho(\PP)\neq\rho(\QQ)$, and we may solve for $t$ and substitute this value along with the given values of $x_i$ and $y_i$ into $g_T$.  After this substitution the $W_j$ corresponding to degenerate triangles vanish. As the $W_i$ and $W_U$ stand for twice the areas, we now divide by $2^{\deg g_T}$ to get the desired integral equation for $u$ over $a_1, \dots, a_n$.
\end{proof}

We conclude this section with a consequence for parallelograms which generalizes a theorem of Monsky.

\begin{cor}\label{cor:sigma/2}
Let $T=\PP\QQ\RR\SS$ be a parallelogram in the Euclidean plane with a dissection into $n$ triangles of areas $a_1, \dots, a_n$.  Let $\sigma$ denote the area of $T$.  
Then $\sigma/2$ is integral over $\Z[a_1,\dots, a_n]$.
\end{cor}

This corollary implies the fact due to Monsky \cite{monsky} that if a square of area $1$ in the Euclidean plane is dissected into $n$ triangles of areas $a_1, \dots, a_n$, then there is a polynomial $f$ with integer coefficients such that $2f(a_1, \dots, a_n) = 1$.  (To see this, take the integral equation for $\sigma/2=1/2$ and multiply by a power of 2 to clear denominators.)  Likewise, Theorem 17 of \cite{chapter3}, which extends Monsky's theorem to handle deformations, can be derived from Theorem 1+.

\section{The Area Map for Trapezoids}

Theorem 1+ tells us that $W_U$ is integral over $\Z[W_1, \dots, W_n]$, i.e., there exists a polynomial $g=g_T\in \Z[U, B_1, \dots, B_n]$, monic in $U$, such that $g(W_U,W_1,\dots,W_n)=0$ in $\thering$.  Assuming that $g$ has been chosen with minimal degree, we will now show that almost all points in the zero set of $g$ are realized as areas of triangles in an actual trapezoidal drawing of $T$.  For this purpose, we introduce a {\it drawing space} $\Trap(T)$ and an {\it area map} for this situation. 

Let $T$ be a combinatorial triangulation of a quadrilateral with corners $\PP\QQ\RR\SS$.  A \emph{drawing} of $T$ is a map $\rho:\Vxs(T)\to \C^2$ that takes $\PP\QQ\RR\SS$ to a trapezoid; this means that the vectors $\QQ-\PP$ and $\RR-\SS$ are linearly dependent.
Let $\Trap=\Trap(T)$ be the space of drawings of $T$.  An open dense subset of $\Trap$ is parameterized by the affine space $X=X(T)$ with coordinates $x_v, y_v$ for all vertices $v$ except $\RR$ and an additional coordinate $t$.
We will keep track of the areas of the triangles of $T$ as well as the area $U$ of the triangle formed by the images of $\PP,\SS,$ and $\QQ$ (even though these vertices probably do not form a triangle of the triangulation); thus let $Y=Y(T)$ denote the projective space with one coordinate for each triangle of $T$ and one additional coordinate $U$.  Now let
$\Area:X \dashrightarrow Y$
be the (rational) area map that records the areas of the triangles in the corresponding coordinates and the area of the triangle $\rho(\PP),\rho(\SS),\rho(\QQ)$ in the $U$ coordinate. 

Let $V=V(T)$ denote the closure of the image of the map $\Area$.  Thus $V\subset Y$ is a rational variety.

\begin{thm}
For any $T$, the variety $V(T)$ is an irreducible hypersurface in $Y$ defined by a homogeneous polynomial $z_T(U, B_1, \dots, B_n)$ that is monic in $U$.
\end{thm}

\begin{proof}
The parameter space $X$ is irreducible, and therefore $V(T)$ is also irreducible. To show $V(T)$ is a hypersurface,
we appeal to the argument from \cite[Theorem 5]{triangles1} that $\Area$ is generically locally injective after modding out by affine transformations.  A dimension count then shows that the image of $\Area$ has codimension 1 in $Y$.

Let $z_T$ be the defining equation of $V(T)$, scaled to have integer coefficients. We wish to show that $z_T$ is monic in $U$.  By Theorem 1+, there exists $g\in \Z[U, B_1\dots, B_n]$ which is monic in $U$ and such that $g(W_U, W_1, \dots, W_n)=0$ in $\mathcal{R} $.  We assume that we have chosen such a $g$ with minimal degree.  Note that $g=g(U, B_1, \dots, B_n)$ vanishes on the image of $\Area$, so $z_T$ divides $g$. 

We now argue that in fact $g=z_T$.
Observe that the $W_i$ are algebraically independent over $\C$, because if there were a dependence $r(W_1, \dots, W_n)=0$, we would have $z_T$ divides $r$, which implies that $z_T$ does not contain the variable $U$.  But then $g$, which is a multiple of $z_T$, would not be monic in $U$, a contradiction. We conclude that $\Z[W_1, \dots, W_n]$ is isomorphic to a polynomial ring, which is a UFD. By Gauss's Lemma, the integral equation $g$ may be chosen to be irreducible as a polynomial in $\Q(W_1, \dots, W_n)[U]$. It follows that $g(U, B_1, \dots, B_n)$ is irreducible in $\Q[U, B_1,\dots, B_n]$. From this we see that  $g=z_T$, and so $z_T$ is monic in $U$, as desired.
\end{proof}

\section{Integrality for Parallelograms}

In this section we prove Theorem \ref{thm2} and Theorem \ref{monicity}. The proofs of these integrality theorems for parallelograms rely on our integrality theorem for trapezoids. 

The polynomial $p_T$ for parallelograms, studied in \cite{triangles1,chapter3,chapter2}, can be linked to the polynomial $z_T$ for trapezoids using a simple geometric observation:  a trapezoid $T=\PP\QQ\RR\SS$ is a parallelogram if  and only if its area is twice the area of triangle $\PP\QQ\SS$.  For a triangulated trapezoid, this condition is represented by the equation $-2U=S$, where $S$ denotes $\sum_{i=1}^n B_i$.  This observation implies the relation
$$
p_T(B_1, \dots, B_n) \ | \ z_T ( -S , 2B_1, \dots, 2B_n ) 
$$
from which we will tease out the monicity of $p_T$.  

To do this, one further fact about $z_T$ is required.

\begin{prop} For any $T$, we have 
$z_T(U,B,0,\dots,0)=\pm U^e(U+B)^f$ for nonnegative integers $e$ and $f$.
\end{prop}

Proving this requires understanding points of $V$ that are not in the image of the area map.  The paper \cite{chapter2} studies this question in a nearly identical context, namely the area map for a triangulated parallelogram.  One main conclusion there is that if $w$ is a point of $V$ then either $w$ is in the image of $\Area$ or else there is a subset of the coordinates that sums nontrivially to $0$.  This conclusion is also valid for the trapezoid area map.

\begin{lem}\label{bubblelemma}
Suppose $w=[u:b_1:\dots :b_n]\in V\setminus\im\Area$. Let $b_0=u$.  Then there is a subset $Z$ of $\{0, \dots, n\}$ such that $\sum_{i\in Z} b_i = 0$, but $b_i\neq 0$ for some $i\in Z$. 
\end{lem}

\begin{proof}
We view $\Area$ as the area map associated to the complex $\hat T=T \cup U$ which is a triangulation of the triangle $\QQ\RR\SS$. The proof is nearly identical to the parallelogram case \cite[Main Theorem 3]{chapter2}. Here are the main points of the argument. We use the language of generating paths and bubbles introduced in \cite[Section 3]{chapter2}.

Suppose $w\in V\setminus\im\Area$. Then there is a generating path for $w$, which is a path $\gamma(s)$ of drawings in $\Trap$ converging to a limiting $\rho\in \Trap$ as $s\rightarrow 0$ and such that $\Area(\gamma(s))\rightarrow w$.

If $\rho$ maps the boundary $\QQ\RR\SS$ to a single point, then $\rho$ contains a bubble.  Otherwise there are two adjacent points of the boundary $V_1$ and $V_2$ such that $\rho(V_1)\neq\rho(V_2)$.  Using an invertible affine transformation we may assume $\rho(V_1)=(0,0)$ and $\rho(V_2)=(1,0)$, and a further affine transformation that converges to the identity as $s\to 0$ fixes $\gamma(s)(V_1)=(0,0)$ and $\gamma(s)(V_2)=(1,0)$. We then rescale vertically so that some vertex is not converging to the $x$-axis.  This produces a new generating path, with a limiting drawing that we still call $\rho$.  By the Elastic Lemma of \cite{chapter2}, $\rho$ must have a bubble.

We conclude that there exists a generating path for $w$ with a bubble. The Bubble Corollary of \cite{chapter2} then asserts that the coordinates inside this bubble sum to zero but are not all zero.  
\end{proof}

We now prove the Proposition.
\begin{proof}  
From Theorem 1+, $z_T$ is monic in $U$ and hence also $z_T(U,B,0,\dots,0)$ is monic in $U$.  Thus it suffices to show that the only zeros of $z_T(U,B,0,\dots,0)$ have $U=0$ or $U=-B$.

Note that $[1:0:\cdots]$ is not in $V$, again since $z_T$ is monic in $U$. So we may assume $B\ne 0$, and suppose $w=[U:1:0:\cdots]\in V$. We will show that $U=0$ or $U=-1$.

If $U= -1$, we are done. Otherwise, by Lemma \ref{bubblelemma}, we have $w\in \im\Area$.  Thus, there is a drawing with $B=1$ and the areas of all other triangles of $T$ equal to $0$. It follows from [2, Corollary 5.6 (1)] that the boundary of $T$ must be drawn as a degenerate trapezoid.  But the vertices of the boundary cannot all be collinear, since then the $B_i$ would sum to $0$. Thus the image of the boundary is a nondegenerate triangle of area $1$, and the four points $\PP, \QQ, \RR, \SS$ map onto the three corners of this triangle.  Thus we see that either the $U$ triangle, $\PP\SS\QQ$, or the $U'$ triangle, $\QQ\SS\RR$, is degenerate.  However $U$ and $U'$ add up to $-\sum B_i$, which equals $-1$.   Hence $U=0$ or $U=-1$.
\end{proof}

We now prove Theorem \ref{monicity} and Theorem \ref{thm2}, in that order. 
\begin{proof}
We first consider the coefficient $\alpha$ of $B_1^d$ in the polynomial 
$$\tilde{z}(B_1, \dots, B_n) = z(-S, 2B_1, \dots, 2B_n).$$
This coefficient $\alpha$ is the same as the coefficient of $B_1^d$ in $z(-B_1, 2B_1,0 \dots, 0)$, which equals $\pm (-B_1)^e B_1^f$ by the Proposition.  Thus $\alpha = \pm 1$. Since $p$ is a factor of $\tilde{z}$, it follows from Gauss's Lemma that $B_1^{d'}$ has coefficient $\pm 1$ in $p$, where $d'$ is the degree of $p$.  This proves Theorem \ref{monicity}.

To prove Theorem \ref{thm2} for triangulations, we view the polynomial $p_T(B_1, \dots, B_n)$ as a polynomial in $B_n$ with coefficients in $\Z[B_1, \dots, B_{n-1}]$.  We have just established that the leading coefficient is $\pm 1$.  Thus $p_T$ provides the required integral equation for $B_n$ over $\Z[B_1, \dots, B_{n-1}]$.

To prove Theorem \ref{thm2} for dissections, apply the poofing argument used in Theorem \ref{thm1} to produce a combinatorial triangulation to which the previous paragraph applies.
\end{proof}

\begin{example*}
The triangulation $T_n$ with vertices $\PP=p_0,p_1,\ldots,p_{n+1}=\RR$, $\QQ$, $\SS$ and triangles $A_i=\SS p_{i-1} p_{i}$ and $B_i=\QQ p_{i} p_{i-1}$ (for $1\le i \le n+1$), called the \emph{diagonal case} in \cite{triangles1}, has 
$$z_{T_n}=\left( \prod\limits_{k=0}^{n+1} \ell_k \right) \left( \frac 1 {\ell_0} - \sum\limits_{k=0}^n \frac {A_{k+1}}{\ell_k \ell_{k+1}} \right)$$
where $\ell_k$ stands for the linear form $A_1+\cdots+A_k+B_1+\cdots+B_k+U$.  Its degree is $n+1$. For example $z_{T_1}=U^2 + 2 U B_1 + U B_2 + U B_4 + B_1^2 + B_1 B_2 + B_1 B_3 + B_1 B_4$.  We then have $z_{T_n}(-S,2A_i,2B_i)=S\cdot p_{T_n}$, where $p_{T_n}$ is computed in \cite{triangles1}.
\end{example*}

\bibliographystyle{unsrt}
\bibliography{IntegralityRefs}

\begin{thebibliography}{1}

\bibitem{triangles1}
Aaron Abrams and James Pommersheim.
\newblock Spaces of polygonal triangulations and {M}onsky polynomials.
\newblock {\em Discrete and Computational Geometry}, 51:132--160, 2014.
\newblock DOI: 10.1007/s00454-013-9553-6.

\bibitem{chapter3}
Aaron Abrams and Jamie Pommersheim.
\newblock Generalized dissections and {M}onsky's theorem.
\newblock {\em Discrete Comput. Geom.}, 67(3):947--983, 2022.

\bibitem{chapter2}
Aaron Abrams and James Pommersheim.
\newblock An illustrated encyclopedia of area relations.
\newblock {\em European {J}ournal of {M}athematics}, 9(49), 2023.

\bibitem{monsky}
Paul Monsky.
\newblock On dividing a square into triangles.
\newblock {\em Amer. Math. Monthly}, 77(2):161--164, 1970.

\bibitem{monskyjepsen}
Charles~H. Jepsen and Paul Monsky.
\newblock Constructing equidissections for certain classes of trapezoids.
\newblock {\em Discret. Math.}, 308(23):5672--5681, 2008.

\bibitem{kasimatisstein}
E.~A. Kasimatis and S.~K. Stein.
\newblock Equidissections of polygons.
\newblock {\em Discrete Math.}, 85(3):281--294, 1990.

\bibitem{atiyahmacdonald}
Michael~Francis Atiyah and I.~G. MacDonald.
\newblock {\em Introduction to commutative algebra.}
\newblock Addison-Wesley-Longman, 1969.

\end{thebibliography}

\end{document}